\newcommand{\CLASSINPUTtoptextmargin}{2cm}
\newcommand{\CLASSINPUTbottomtextmargin}{2cm}
\documentclass[letterpaper, 10 pt, journal, twoside]{ieeetran}
\usepackage{cite}
\usepackage{amsmath,amssymb,amsfonts,bm}
\usepackage{graphicx}
\usepackage{textcomp}
\usepackage{xcolor}
\usepackage{color}
\usepackage{amsthm}
\usepackage{amsmath}
\usepackage[labelsep=period]{caption}
\usepackage{graphics} 
\usepackage{epsfig}
\usepackage{times} 
\usepackage{textcomp}
\usepackage{multicol}
\usepackage{float}
\usepackage{hyperref}
\usepackage{multirow}
\usepackage{subcaption}
\usepackage{graphicx}
\usepackage[top=54pt, left=54pt, right=54pt, bottom=54pt]{geometry}

\newtheorem{thm}{Theorem}
\newtheorem{prop}[thm]{Proposition}

\newtheorem{rem}{Remark}
\newtheorem{lem}[thm]{Lemma}

\usepackage{enumerate}

\usepackage{algorithm}
\usepackage{algorithmicx}
\usepackage{algpseudocode}
\usepackage{flushend}

\DeclareMathOperator*{\argmin}{arg\,min}

\usepackage[flushleft]{threeparttable}

\title{\LARGE \bfseries Minimax Linear Optimal Control of Positive Systems}

\author{Alba Gurpegui, Emma Tegling and Anders Rantzer
\thanks{The authors A. Gurpegui, E. Tegling and A. Rantzer are with the Department of Automatic Control and the ELLIIT Strategic Research Area at Lund University, Lund, Sweden. Email: \href{mailto:alba.gurpegui@control.lth.se}{alba.gurpegui@control.lth.se}, \href{mailto:emma.tegling@control.lth.se}{emma.tegling@control.lth.se}, 
\href{mailto:anders.rantzer@control.lth.se}{anders.rantzer@control.lth.se}. This work is partially funded by the Wallenberg AI, Autonomous Systems and Software Program (WASP) and the European Research Council (ERC) under the European Union's Horizon 2020 research and innovation programme under grant agreement No 834142  (ScalableControl).}%

}


\begin{document}
\IEEEoverridecommandlockouts
\IEEEpubid{\makebox[\columnwidth]{XXX-X-XXXX-XXXX-X/XX/\$XX.XX~\copyright{}2024 IEEE \hfill} \hspace{\columnsep}\makebox[\columnwidth]{ }}

\maketitle
\thispagestyle{empty} 
\pagestyle{empty}    


\maketitle
\begin{abstract}
We present a novel class of minimax optimal control problems with positive dynamics, linear objective function and homogeneous constraints. The proposed problem class can be analyzed with dynamic programming and an explicit solution to the Bellman equation can be obtained, revealing that the optimal control policy (among all possible policies) is linear. This policy can in turn be computed through standard value iterations. 
Moreover, the feedback matrix of the optimal controller inherits the sparsity structure from the constraint matrix of the problem statement. This permits structural controller constraints in the problem design and simplifies the application to large-scale systems. We use a simple example of voltage control in an electric network to illustrate the problem setup. 
\end{abstract}
\section{Introduction}
\subsection{Motivation}
Optimal control problems with minimax objectives are ubiquitous in control theory and engineering. They provide a powerful framework for modeling and solving problems that involve competition and uncertainty~\cite{Tamer_Basar, rockafellar}. These problems appear in contexts such as robust control, game theory, and multi-agent systems. In our particular case, they are used to design control systems that are robust to uncertainties and disturbances, with the objective to minimize the worst-case performance of the system. Finding solutions to these types of problems can be a challenging task, especially when dealing with large-scale systems. 

In this paper, we present a novel class of minimax optimal control problems, with positive dynamics, linear objective function, and linear, homogeneous, constraints. Inspired by the explicit solution presented in~\cite{mainArt} for the minimization case, the solution of this class of problems is based on dynamic programming, with an explicit solution to the problem's Bellman equation. Thus, it is possible to find the optimal control policy, which, for this class of problems, is manifested as a linear feedback policy that minimizes 
the objective function over the system's trajectory when subjected to the worst-case disturbance or uncertainty that is homogeneous in the system state. 

To understand the relevance of our problem class, recall that a linear system is positive if the state and output remain nonnegative as long as the initial state and the inputs are nonnegative. This type of dynamics have gained attention in control theory literature because of all the technological and physical phenomena that can be captured 
by positive dynamics. Classical books on the topic are~\cite{BermanBook} and~\cite{Luenberger}. In the latter, David Luenberger in 1979 devotes a chapter to positive dynamical systems, which is considered by many the initiation of “Positive System Theory.”  Significant research has been conducted to represent natural extensions of the class of positive systems, for instance positive systems with delays~\cite{Eibhara_steady_state_delay_interc_positive}, positive switched systems~\cite{Blanchini_switched_linear_positive} and monotone systems~\cite{Smith1995MonotoneDS}. Positive systems theory has been useful to describe dynamical systems in a wide range of applications, such as, biology, ecology, physiology and pharmacology~\cite{17tutorial, 19tutorial, 37tutorial, 38tutorial,39tutorial,44tutorial}, thermodynamics~\cite{12tutorial, 38tutorial}, epidemiology~\cite{1tutorial, 40tutorial,53tutorial}, econometrics~\cite{55tutorial}, filtering and charge routing networks or power systems~\cite{6tutorial,7tutorial,10tutorial}. 

 One of the main advantages of positive systems is that stability can be verified using linear Lyapunov functions~\cite{Blanchini_lyapunov}, making this
class of systems more tractable in a large scale setting because of their computational scalability~\cite{tutorial_anders_valcher,EbiharaPeucelleArzelierTAC2017}. Another well known advantage is that linear controllers $u=-Kx$ can be designed with sparsity constraints on $K$. See for example \cite{Briat2013, Tanaka+11}. This paper is optimizing sparse controllers of the same form, but with one important difference compared to past literature: \emph{Dynamic programming is carried out without a priori constraints on linearity or sparsity. Instead these properties are a consequence of the optimization critera and constraints.} Hence it is possible to conclude that no nonlinear nonsparse controller can ever achieve a lower value of the cost. The minimax optimization is also related to past work on gain minimization \cite{Ebihara1,Briat2013}, but again the dynamic programming approach is different. More powerful conclusions are obtained at the expense of more restrictive assumptions.

\subsection{Problem Setup}
We present the optimal control problem of this paper as a discrete-time, infinite-horizon, minimax optimal control problem with nonnegative cost and positive dynamics with linear, continuous objective function and  homogeneous constraints,  
\begin{align}\label{optprob}
    \underset{\mu}{\mathrm{inf}}\hspace{1mm}\underset{w}{\mathrm{\max}}&\sum_{t = 0}^{\infty} \left[s^{\top}x(t)+r^{\top}u(t)-\gamma^{\top}w(t) \right ]\\
    \mathrm{subject}& \hspace{1mm} \mathrm{to} \notag \\
    &x(t+1) = Ax(t)+Bu(t)+Fw(t), \notag \\
    &u(t)=\mu(x(t)) \hspace{1mm}; \hspace{2mm} x(0)=x_{0} \notag \\
    &\left | u \right |\leq Ex \hspace{1mm}; \hspace{1mm} \left | w \right |\leq Gx \notag
\end{align}
where $x$ represents the $n$-dimensional vector of state variables, $u$ the $m$-dimensional control variable, $w$ the $l$-dimensional disturbance, $\mu$ is any, potentially nonlinear, control policy, $E$ prescribes the structure of the control action and $G$ is assumed to determine the linear dependency of the disturbance and the state. The objective is to minimize the worst-case cost over all possible control strategies. 

\subsection{Notation} Let $\mathbb{R}_{+}$ denote the set of nonnegative real numbers. The inequality $X > Y$ $(X \geq Y)$ mean that all the elements of the matrix $(X-Y)$ are positive (nonnegative).  A matrix $X$ is called positive if all the elements of $X$ are nonnegative but at least one element is nonzero. The notation $\left | X \right |$ means elementwise absolute value.
\section{MAIN RESULT}
In this section we state and prove the main theorem of the paper, which gives an explicit solution to the presented class of linear minimax optimal control problems \eqref{optprob}.
\begin{thm}\label{MainThm} Let $A \in \mathbb{R}^{n\times n}$, $B = \left [ B_{1}^{\top}, \ldots, B_{m}^{\top} \right ]^{\top} \in \mathbb{R}_{+}^{m \times n}$, $F \in \mathbb{R}^{n \times l}$, $E = \left [ E_{1}^{\top}, \ldots, E_{m}^{\top} \right ]^{\top} \in \mathbb{R}_{+}^{m \times n}$,  $G \in \mathbb{R}_{+}^{l \times n}$, $s \in \mathbb{R}^{n} $, $r \in \mathbb{R}^{m} $, $\gamma \in \mathbb{R}^{l} $. Suppose that 
    \begin{align}\label{ass1}
        A &\geqslant \left | B \right |E + \left | F \right |G \\
        s&\geqslant  E^{\top} \left | r \right | - G^{\top} \left | \gamma \right |.
        \label{ass2}
    \end{align}
Then the following statements are equivalent:
\begin{enumerate}[(i)]
    \item The optimal control problem~\eqref{optprob}, has a finite value for every $x_{0} \in \mathbb{R}_{+}^{n}$.
    \vspace{1mm}
    
    \item The recursive sequence $\left \{ p_{k} \right \}_{k=0}^{\infty}$ with $p_{0}=0$ and
    \begin{multline} \label{sequence_main_p}
     p_{k} = s + A^{\top}p_{k-1}  -E^{\top} \left | r+B^{\top}p_{k-1} \right | \\ + G^{\top}\left | -\gamma + F^{\top}p_{k-1} \right |
    \end{multline}
has a finite limit.
    \vspace{1mm}
    
\item There exists $p\in \mathbb{R}_{+}^{n}$ such that
{\begin{align}\label{p_def}
    p\! = \!s + A^{\top}p - E^{\top}\!\left | r \!+\! B^{\top}p \right | + G^{\top}\!\left | -\gamma \!+ \!F^{\top}p \right |.
\end{align}}
\end{enumerate}
\normalsize{If (iii) is true then~\eqref{optprob} has the minimal, finite, nonnegative value $p^{T}x_{0}$, with $p$ being the limit of the recursive sequence $\left \{ p_{k} \right \}_{k=0}^{\infty}$ in $(ii)$. Moreover, the control law $u(t) = -Kx(t)$, is optimal when }
\begin{align} \label{Kform}
    K := \begin{bmatrix}
\mathrm{sign}(r_{1}+ p^{\top}B_{1})E_{1}\\ 
\vdots \\ 
\mathrm{sign}(r_{m}+p^{\top}B_{m})E_{m}
\end{bmatrix}.
\end{align}
\end{thm}
\vspace{2mm}
\begin{rem}
    Even though the optimization is done without any \emph{a priori} assumption of linearity of the control policy, the resulting optimal controller is linear. 
\end{rem}
\begin{rem}
The condition~\eqref{ass1} ensures the invariance of the positive orthant under the system dynamics.  
The second condition~\eqref{ass2} will be needed when applying Proposition~\ref{propAux} in the Appendix to our objective function $g(x,u,w) = s^{T}x + r^{T}u -\gamma ^{T}w$.
\end{rem}
\begin{rem}
In \eqref{Kform}, it can be observed that the sparsity structure of the control gain $K$ is directly determined by the $E$ matrix. The sparsity of $E$ is in turn determined by the problem designer and may capture limitations in actuation and sensing. 
\end{rem} 
\begin{proof}
    The general nonlinear minimax optimal control problem~\eqref{optprob_g} presented in the Appendix~\ref{APPENDIX} reduces to our problem set up~\eqref{optprob} if 
\begin{align*}
    f(x,u,w) &:= Ax+Bu+Fw \\
    g(x,u,w)&:= s^{\top}x+r^{\top}u - \gamma^{\top}w.
\end{align*}
Furthermore, under condition~\eqref{ass2} we observe that
\begin{align*}
    &\underset{\left | w \right |\leq Gx}{\max}\left[s^{\top}x + r^{\top}u- \gamma^{\top}w \right]\\
    &= s^{\top}x + r^{\top}u+ \left | \gamma \right |^{\top}Gx 
    \\
    &\geq\left( E^{\top} \left | r \right | - G^{\top} \left | \gamma \right |\right)^{\top}x+r^{\top}u+ \left | \gamma \right |^{\top}Gx\\
    &\geq  - \left( G^{\top} \left | \gamma \right | \right)^{\top}x +\left | \gamma \right |^{\top}Gx  = 0.
\end{align*}
Therefore, 
\begin{align*}
    \small{\underset{\left | w \right |\leq Gx}{\max}\left[g(x,u,w) \right]=\underset{\left | w \right |\leq Gx}{\max}\left[s^{\top}x + r^{\top}u- \gamma^{\top}w \right] \geq 0} 
\end{align*}
\normalsize{as} required in Lemma~\ref{LEMAA} in the Appendix~\ref{APPENDIX}. Now, it is clear that $(i)$ in Theorem~\ref{MainThm} is equivalent to $(i)$ in Lemma~\ref{LEMAA} in the Appendix~\ref{APPENDIX}. Next, we will verify that the recursive sequence in $(ii)$ is equivalent to $(ii)$ in Lemma~\ref{LEMAA}. To prove this we use induction over $p^{\top}_{k}x=J_{k}(x)$. By definition, it is direct that $p_{0}^{\top}x =0=J_{0}(x)$ for all $x$. For the induction step we assume that $p^{\top}_{k}x=J_{k}(x)$. Now we want to prove that $p^{\top}_{k+1}x=J_{k+1}(x)$. From~\eqref{sequence_J} and the induction hypothesis we have that 
\begin{align*}
    &J_{k+1}(x)=\underset{u }{   \mathrm{\min} }\hspace{1mm}\underset{w }{\mathrm{\max}}\left [ g(x,u,w)+J_{k}(f(x,u,w)) \right ]\\
    &=\underset{u}{\mathrm{min}} \hspace{1mm} \underset{w}{\mathrm{max}}\left [ s^{\top}x + r^{\top}u-\gamma^{\top}w+ p^{\top}_{k}(Ax+Bu+Fw))\right ]\\
    &=s^{\top}x + p_{k}^{\top}Ax+\underset{\left | u \right |\leq Ex}{\mathrm{min}}\left [ r^{\top}u+ p^{\top}_{k}Bu \right ] \\ 
    &~~~~~~~~~~~~~~~~~~~~~~~~~~~~~~~~~~~~+ \underset{\left | w \right |\leq Gx}{\mathrm{max}}\left [-\gamma^{\top}w+ p^{\top}_{k}Fw \right ]\\
    &=s^{\top}x + p_{k}^{\top}Ax  - \left | r+B^{\top}p_{k} \right |^{\top}Ex + \left | -\gamma + F^{\top}p_{k} \right |^{\top}Gx\\
    &=p^{\top}_{k+1}x.
\end{align*}
\noindent Therefore, $p_{k}^{\top}x=J_{k}(x)$ for all $k \in \mathbb{N}$ and all $x \in \mathbb{R}^{n}_{+}$, and $p^{\top}x=J^{*}(x)$ for all $x$. Hence, $(ii)$ and $(iii)$ both in Theorem~\ref{MainThm} and Lemma~\ref{LEMAA} are equivalent. Furthermore, note that under homogeneous constraints the linearity of $J_{k}$ is preserved during value iteration.

Because $(i)$, $(ii)$ and~$(iii)$ in this theorem and in Lemma~\ref{LEMAA} are equivalent, the proof of equivalence between $(i)$, $(ii)$ and $(iii)$ in Theorem~\ref{MainThm} follows from the proof of Lemma~\ref{LEMAA} in the Appendix~\ref{APPENDIX}.

To finish this proof it is just left to give an expression for the optimal control policy $u(t)=\mu(x(t))$, 
\begin{align*}
    \mu(x) &= \argmin_{ | u|\leq Ex} 
    \left [ s^{\top}x+r^{\top}u-\gamma^{\top}w + p^{\top}(Ax+Bu+Fw) \right ]\\
   &= \argmin_{ | u|\leq Ex} \sum_{i=1}^{m}\left[ \left ( r_{i}+p^{\top}B_{i} \right )u_{i} \right].
\end{align*}
\\
Finally, since for all $i = 1 \hspace{1mm}... \hspace{1mm} m$ the inequality $\left|u \right | \leq Ex$ restricts $u_{i}$ to the interval $\left [ -E_{i}x, E_{i}x \right ]$, the minimum is attained when $(r_{i}+p^{\top}B_{i})$ and $u_{i}$ have opposite signs. Thus, $u_{i} = -\mathrm{sign}(r_{i}+p^{\top}B_{i})E_{i}$ for all $i = 1 \hspace{1mm}... \hspace{1mm} m$.
\end{proof}
Next, we present an example of a simple scenario where the main theorem is applied to a network problem -- a DC power network -- to illustrate the role of the asssumptions and contraints. 
\begin{figure}
    \includegraphics[scale=0.21]{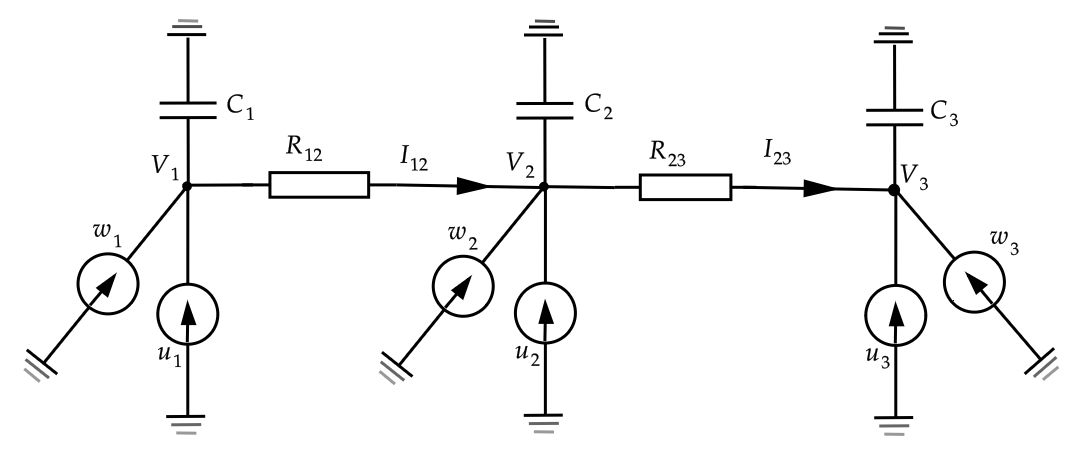}
    \caption{ Example of a DC network consisting of 3 terminals (buses) and 5 lines. The controls $u_i$ are used to control the voltage when the system is subjected to the disturbances $w_i$.}
    \label{emmadc}
\end{figure}
\section{Example: Optimal voltage control in a DC power network}
The optimal control problem~\eqref{optprob} admits sparsity constraints on the controller, making it particularly useful for problems defined over network graphs. Here, we 
consider a simple example of voltage control in a 
DC (i.e., direct current) power network. 
Here, the nodes represent voltage source converters with positive voltage dynamics, interconnected through resistive lines. The model can, for example, capture an envisioned multi-terminal high-voltage DC network, whose design aims to transmit power over long distances while maintaining low losses~\cite{VANHERTEM2010,EmmaMTDC} or a simplified DC distribution network~\cite{Karlsson2003DC}. The (continuous) voltage dynamics at the DC bus~$i$ (node~$i$) is given by: 
{\begin{align}\label{MTDC_EQ} \nonumber
    C_{i}\dot{V}_{i}(t) &= -\sum_{j=1}^{n}I_{ij}+u_{i}(t)+w_{i}(t) \\ &= - \sum_{j=1}^{n} \frac{1}{R_{ij}}(V_{i}(t)-V_{j}(t)) + u_{i}(t) + w_{i}(t) ,
\end{align}}
for all $i = 1,2,\ldots,n$, where $u_{i}$ denotes the controlled injected current, $R_{ij}$ the resistance of transmission line~$(i,j)$ (with $R_{i,j}=\infty$ if there exists no line connecting nodes $i$ and $j$),  and $C_{i}$ is the total capacitance at bus~$i$.\footnote{Any line capacitances can for the purpose of this example be absorbed in to the buses. }  
We have also included the disturbance current $w_i$, arising from variations in local generation and load.  
Defining the vector $V=\left[V_{1}(t),...,V_{n}(t) \right]^{\top}$, $u $ and $w$ analogously, and $C = \mathrm{diag}(\left[C_{1},...C_{n}\right])$, we may write~\eqref{MTDC_EQ} on vector form as  
\begin{align}\label{MTDCEQ2}
    C\dot{V}(t)= -\mathcal{L}_{R}V(t)+ u(t) +w(t).
\end{align} 
Here, $\mathcal{L}_{R}$ is the weighted Laplacian matrix of the graph representing the transmission lines, whose edge weights are given by the conductances $\frac{1}{R_{ij}}$, i.e., 
\begin{align*}
    \left [ \mathcal{L_{R}} \right ]_{i,j} = \left\{\begin{matrix}
-\frac{1}{R_{i,j}} \hspace{2mm} \mathrm{if} \hspace{2mm} i \neq j\\ 
\sum_{j=1}^{n} \frac{1}{R_{i,j}} \hspace{2mm} \mathrm{if} \hspace{2mm} i = j
\end{matrix}\right. .
\end{align*}
Note that $-\mathcal{L}_{R}$ is Metzler, and the system~\eqref{MTDCEQ2} thus positive.
\\
The dynamics in~\ref{MTDCEQ2} can be discretized as $$C ( {V(\tau h+h)-V(\tau h)}  )/h=-\!\mathcal{L_{R}}V(\tau h)+u(\tau h)+w(\tau h).$$
Setting $t = \tau h$ and re-defining the state $x(t) = V(\tau h)$ gives the discrete-time dynamics 
\begin{multline} \label{discrete_dyn}
     x(t+1) = \left [ I-hC^{-1} \mathcal{L}_{R} \right ]x(t) + hC^{-1}u(t) \\+ hC^{-1}w(t).
\end{multline}
Now, we formulate the optimal control problem~\eqref{optprob} for the dymamics~\eqref{discrete_dyn}: 
\begin{align}\label{optprobMTDC}
    \underset{\mu}{\mathrm{inf}}\hspace{1mm}\underset{w}{\mathrm{\max}}&\sum_{t = 0}^{\infty} \left[s^\top{}x(t)+r^{\top}u(t)-\gamma^{\top}w(t)  \right ]\\
    \mathrm{Subject} &\hspace{1mm} \mathrm{to} \notag \\
    &x(t+1) = \left [ I-hC^{-1} \mathcal{L}_{R} \right ]x(t) \\ \notag &~~~~~~~~~~~~~~~~+hC^{-1}u(t) + hC^{-1}w(t) \notag\\
    &u(t)=\mu(x(t)) \hspace{1mm}; \hspace{2mm} x(0)=x_{0} \notag \\
    &\left | u \right |\leq Ex \hspace{1mm}; \hspace{1mm} \left | w \right |\leq Gx \notag
\end{align}
Identifying $A$ and $B$, condition (\ref{ass1}) reads
\begin{align}\label{ass1_ex}
        \left [ I-hC^{-1} \cdot \mathcal{L}_{R} \right ]&\geqslant hC^{-1}E + hC^{-1}G. 
\end{align}
Clearly, the right hand side of the inequality must inherit the zero pattern of $\mathcal{L}_{R}$, i.e. its sparsity pattern. In other words, the disturbances and control signal must be compatible 
with the physical network structure and depend only on connected nodes. $E$, or $G$ can, however be more sparse than $\mathcal{L}_{R}$.
Furthermore,~\eqref{ass1_ex} reveals that the diagonals of the left hand side must satisfy 
\begin{align*}
(e_{ii}+g_{ii})+\sum_{j=1}^{N}\frac{1}{R_{ij}}\leq \frac{C_{i}}{h}
\end{align*}
This can always be satisfied by making $h$ sufficiently small. However, the off-diagonals reveal conditions on $e_{ij}$, $g_{ij}$ that depend on the line resistances $R_{ij}$ in a manner best illustrated by~\eqref{exeq}. 

In Fig.~\ref{emmadc} a 3-terminal DC power network system is introduced. For this network, the condition~\eqref{ass1_ex} reads: 
\begin{small}
\small{\begin{align} \label{exeq}
        &\begin{bmatrix}
        1-\sum_{j=1}^{3}\frac{h}{R_{1,j}\cdot C_{1}} & \frac{h}{R_{1,2}\cdot C_{1}}  & 0  \\ 
       \frac{h}{R_{2,1}\cdot C_{2}} & 1-\sum_{j=1}^{3}\frac{h}{R_{2,j}\cdot C_{2}}  & \frac{h}{R_{2,3}\cdot C_{2}} \\ 
       0 &\frac{h}{R_{3,2}\cdot C_{3}}  & 1- \sum_{j=1}^{3}\frac{h}{R_{3,j}\cdot C_{3}} \end{bmatrix} \notag \\
       &\geqslant \begin{bmatrix}
        \frac{h}{C_{1}}(e_{1,1} + g_{1,1}) & \frac{h}{C_{1}}(e_{1,2} + g_{1,2}) & 0 \\ 
        \frac{h}{C_{2}}(e_{2,1} + g_{2,1}) & \frac{h}{C_{2}}(e_{2,2} + g_{2,2}) &\frac{h}{C_{2}}(e_{2,3} +g_{2,3})\\ 
        0 & \frac{h}{C_{3}}(e_{3,2} + g_{3,2}) & \frac{h}{C_{3}}(e_{3,3} + g_{3,3})
        \end{bmatrix} .
\end{align}}
\end{small}
\normalsize{This element-wise matrix inequality shows necessary constraints on the elements of $E$ and $G$.}
\newline
\\
In parallel, condition~\eqref{ass2} means that $E$ and $G$ need to satisfy 
\begin{align*}
    s\geqslant  E^{\top} \left | r \right | - G^{\top} \left | \gamma \right |.
\end{align*}
Here, the structure of $E$ determines the states available to the local current controllers and $G$ the manner in which disturbances enter the system.
\newline
\\
Particularly, in our 3 terminal DC power network we need the problem design to satisfy
\begin{small}
\begin{align}
\begin{bmatrix}
s_{1}\\ 
s_{2}\\ 
s_{3}
\end{bmatrix}& \! \geqslant \!  \begin{bmatrix}
e_{1,1} &  e_{2,1}& 0\\ 
e_{1,2} & e_{2,2} & e_{3,2}\\ 
0 & e_{2,3} & e_{3,3}
\end{bmatrix} \! \left|\begin{bmatrix}
r_{1}\\ 
r_{2}\\ 
r_{3}
\end{bmatrix}\right|\! -\!\begin{bmatrix}
g_{1,1} &  g_{2,1}& 0\\ 
g_{1,2} & g_{2,2} & g_{3,2}\\ 
0 & g_{2,3} & g_{3,3}
\end{bmatrix} \! \left|\begin{bmatrix}
\gamma_{1}\\ 
\gamma_{2}\\ 
\gamma_{3}
\end{bmatrix} \right|. \notag
\end{align}
\end{small}
\newline
In this example, the resulting optimal controller~\eqref{Kform} can take  8 different configurations depending on the sign of each element of the parameter $r$. If $r$ is positive, because in this problem $B=hC^{-1}$ is positive, all the signs of the rows in $K$ are positive so that $u(t) = -Ex(t)$ becomes optimal.  However, if $r$ is not positive, it is possible to use this parameter to modify the signs of the rows in the resulting control action.
\section{CONCLUSIONS}
In this paper, we have extended the optimal control problem class presented in \cite{mainArt} by exploring the minimax worst-case. Specifically, we have derived a solution $p$ for this novel class of optimal control problems using value iteration. Our resulting optimal controller bears resemblance to the one derived in~\cite{BLANCHINI_REV2} for a continuous-time problem, however, the homogenous constraints we impose allow for a prescribed controller structure.  
We believe that our approach is an interesting first step in analyzing a new class of optimal control problems, but that there is still room for improving our methodology. Considering the high computational cost of value iteration, the exploration of alternatives such as policy iteration \cite{BertsekasVIPI}, offers a compelling direction for future research. 

Our results demonstrate that this class of problems can be scaled to large dynamical systems efficiently, partly since the sparsity of the optimal feedback is directly related to the constraints imposed in the problem statement. We demonstrated the problem setup on a simple example of a 3-node power network, but the same method can be applied when the network scales. 

The explicit solution we obtain to the Bellman equation appears to rely on the homogeneity of the constraints on $u$ and $w$. We believe, however, that constant upper bounds on the signals can be accounted for, and this is the subject of ongoing research.  
\section*{Acknowledgments}
\addcontentsline{toc}{section}{Acknowledgments}
We would like to thank Richard Pates and anonymous reviewers for their insightful comments.
\bibliographystyle{unsrt} 

\section{Appendix}\label{APPENDIX}
\subsection{General Problem Set Up} 
\normalsize{We define, a general discrete time, infinite horizon, minimax optimal control problem with continuous cost function and constraints as}
\begin{align}\label{optprob_g}
    \underset{\mu}{\mathrm{inf}}\hspace{1mm}\underset{w}{\mathrm{\max}}&\sum_{t = 0}^{\infty} g(x(t), u(t), w(t)) \\
    \mathrm{subject}& \hspace{1mm} \mathrm{to} \notag \\
    &x(t+1) = f(x(t), u(t), w(t)), \notag \\
    & x(t)\in X; \hspace{2mm} x(0)=x_{0}\hspace{1mm}; \hspace{2mm} u(t)=\mu(x(t)) \notag \\
    &u(t)\in U(x(t)); \hspace{2mm} w(t)\in W(x(t)) \notag
\end{align}
where $f: X \times U \times W \rightarrow X$, $x$ represents the vector of $n$-dimensional state variables, $u$ the $m$-dimensional control variable and $w$ the $q$-dimensional disturbance. 
\begin{lem} \label{LEMAA}
Suppose 
\begin{align*}
    \underset{w\in W(x)}{\mathrm{max}}\left[g(x,u,w) \right] \geq 0
\end{align*}
$\forall x \in X$, $\forall u \in U(x)$. \normalsize{Then, the following statements are equivalent. }
\begin{enumerate}[(i)]
    \item  The general optimal control problem in~\eqref{optprob_g} has a finite value for every $x_{0} \in \mathbb{R}_{+}^{n}$. 
    \vspace{1mm}
    
    \item The recursive sequence $\left \{ J_{k} \right \}_{k=0}^{\infty}$ with $J_{0} =0$ and
    \begin{align}\label{sequence_J}
         J_{k}(x)=\underset{u }{\mathrm{\min} }\hspace{1mm}\underset{w }{\mathrm{\max}}\left [ g(x,u,w)+J_{k-1}(f(x,u,w)) \right ]
    \end{align}
    has a finite limit $\forall x \in X$.
    \vspace{1mm}
    
\item The Bellman equation \begin{align}\label{Bellman_EQ_g}
    J^{*}(x)&=\underset{u}{\mathrm{min}} \hspace{1mm} \underset{w }{\mathrm{max}}\left [ g(x,u,w) + J^{*}(f(x,u,w))\right ]
\end{align} 
has nonnegative solution $J^{*}(x)$, $\forall x \in X$.
\end{enumerate}
\end{lem}
\begin{proof}
    Note that, by dynamic programming, the increasing monotone recursive sequence $\left \{ J_{k} \right \}_{k=1}^{\infty}$ (see Proposition \ref{propAux}) defined in (\ref{sequence_J}) also satisfies 
\begin{align}\label{J_k}
    J_{k}(x) &= \underset{\mu}{\mathrm{inf}}\max_{w}\sum_{t = 0}^{k} \left[ g(x(t), u(t), w(t))\right]. 
\end{align}
\normalsize{To prove the equivalence we will prove implications $(i)\Rightarrow(ii)$, $(ii)\Rightarrow(iii)$, $(iii)\Rightarrow(ii)$, $(iii)\Rightarrow(i)$.

$(i)\Rightarrow(ii)$ Assume $(i)$.  Then, in~\eqref{J_k} when $k \rightarrow \infty$ we get $J_{k}(x) < \alpha < \infty$} for all $x$, with $\alpha$ representing a finite value for~\eqref{optprob_g}. Hence, the recursive sequence~\eqref{sequence_J} has a finite limit.

$(ii)\Rightarrow(iii)$ Assume $(ii)$. Taking  the limit when $k \rightarrow \infty$ on both sides of the equation~\eqref{sequence_J} we get $(iii)$.

$(iii)\Rightarrow(ii)$ Assume $(iii)$. We want to prove that $\lim_{k \rightarrow \infty}J_{k}(x) < \infty$ for all $x$ with $J_{k}(x)$ defined in~\eqref{J_k}. To achieve this, we use induction over $J_{k}(x) \leq J^{*}(x) < \infty$. It is clear that $0=J_{0}(x) \leq J^{*}(x)$ with $J^{*}$ nonnegative by definition. For the induction step we assume that $J_{k}(x) \leq J^{*}(x)$. We want to prove that $J_{k+1}(x) \leq J^{*}(x)$. From the induction hypothesis it is direct that
\begin{align*}
    J_{k+1}(x)&=\min_{u} \max_{w}\left [ g(x,u,w)+J_{k}(f(x,u,w)) \right ]\\
    & \leq \min_{u} \max_{w}\left [ g(x,u,w)+J^{*}(f(x,u,w)) \right ]\\
    &= J^{*}(x).
\end{align*}
Thus, $J_{k+1}(x) \leq J^{*}(x)$ for all $x$, and $\lim_{k \rightarrow \infty}J_{k}(x) <\infty$ for all $k$ and for all $x$, as we wanted to prove.

$(iii) \Rightarrow (i)$ Assume $(iii)$. Define for all $x$
\begin{align*} 
\mu^{*}(x) = \argmin_{u} \max_{w}\left \{ g(x,u,w)+J^{*}(f(x,u,w)) \right \}
\end{align*}
such that 
\begin{align*} 
J^{*}(x) &=\min_{u} \max_{w}\left [ g(x,u,w) + J^{*}(f(x,u,w))\right ].
\end{align*}
Indeed, from~\eqref{J_k} and implication $(iii) \Rightarrow (ii)$ it can be observed that
\begin{align*} 
\underset{w}{\mathrm{max}}\sum_{t=0}^{k}\left[g(x, \mu^{*}(x),w) \right] &\leq \underset{\mu}{\mathrm{inf}}\max_{w}\sum_{t = 0}^{k} \left[ g(x(t), u(t), w(t))\right]\\
&= J_{k}(x) \leq J^{*}(x) < \infty.
\end{align*}
Hence,
\begin{align*} 
\underset{w}{\mathrm{max}}\sum_{t=0}^{k}\left[g(x, \mu^{*}(x),w)\right] \leq J^{*}(x)
\end{align*} 
for all $k$ and for all $x$. This proves $(i)$.
\end{proof}
\begin{prop}\label{propAux} Let 
\begin{align*}
    \underset{w}{\mathrm{max}}\left[g(x, u, w)\right] \geq 0
\end{align*}
$\forall x \in X$ and $\forall u \in U(x)$. Then, the recursive sequence $\left \{ J_{k}(x) \right \}_{k=0}^{\infty}$ with $J_{0}(x)=0$ and
\begin{align}
    J_{k}(x) =  \min_{u} \max_{w} \left [ g(x, u, w)+J_{k-1}(f(x,u,w)) \right ]
\end{align}
satisfies that  $0\leq J_{0}(x)\leq J_{1}(x)\leq J_{2}(x)\leq ...$ for all $ x \in X$ and all $k \in \mathbb{N}$.
\end{prop}
\begin{proof}
To prove this we use induction over $J_{k}(x)$. From the proposition statement $J_{0}(x)=0$ gives
\begin{align*}
    J_{1}(x) &= \min_{u} \max_{w}\left [ g(x, u, w)+J_{0}(f(x, u, w)) \right ]=\\
    &=  \min_{u} \max_{w}\left [ g(x, u, w)\right ] \geq 0 = J_{0}(x)
\end{align*}
for all $x$. For the induction step we assume that $J_{k}(x) \geq J_{k-1}(x)$. We then want to prove that $J_{k+1} \geq J_{k}$. From \eqref{sequence_J} and the induction hypothesis we have that
\begin{align*}
    J_{k+1}(x) &= \min_{u \in U(x)} \max_{w \in W(x)}\left [ g(x, u, w)+J_{k}(f(x, u, w)) \right ]\\
    &\geq \min_{u \in U(x)} \max_{w \in W(x)}\left [ g(x, u, w)+J_{k-1}(f(x, u, w)) \right ]\\
    &=J_{k}(x) 
\end{align*}
Thus, $J_{k+1}(x)\leq J_{k}(x)$ for all $k$ and for all $x$. 
\end{proof}
\end{document}